\documentclass{ifacconf}

\usepackage{graphicx}      % include this line if your document contains figures
\usepackage{natbib}        % required for bibliography

\usepackage{amsmath} % assumes amsmath package installed
\usepackage{amssymb}  % assumes amsmath package installed

\usepackage{lineno}
\usepackage{xcolor}

\DeclareMathOperator{\dom}{\mathtt{dom}}
\DeclareMathOperator{\interior}{\mathtt{int}}

\DeclareMathOperator{\rge}{\mathtt{rge}}
\DeclareMathOperator{\A}{\mathcal{A}}
\let\H\relax
\DeclareMathOperator{\H}{\mathcal{H}}
\DeclareMathOperator{\K}{\mathcal{K}}
\DeclareMathOperator{\PD}{\mathcal{P}\mathcal{D}}
\DeclareMathOperator{\N}{\mathbb{N}}
\DeclareMathOperator{\R}{\mathbb{R}}

\newtheorem{definition}{Definition}[section]

\newtheorem{proposition}{Proposition}[section]
\newtheorem{example}{Example}[section]
\newtheorem{lemma}{Lemma}[section]
\newtheorem{remark}{Remark}[section]

\begin{document}
\begin{frontmatter}

\title{Prolongation and stability of Zeno solutions to hybrid dynamical systems\thanksref{footnoteinfo}}

\thanks[footnoteinfo]{This work was supported by the German Federal Ministry
 of Education and Research (BMBF) as a part of the research project ''LadeRamProdukt''.}

\author[First]{Sergey Dashkovskiy}
\author[Second]{Petro Feketa}

\address[First]{University of W{\"u}rzburg,
   W{\"u}rzburg, Germany \\ (e-mail: sergey.dashkovskiy@mathematik.uni-wuerzburg.de)}
\address[Second]{University of Applied Sciences Erfurt,
   Erfurt, Germany \\ (e-mail: petro.feketa@fh-erfurt.de)}

\begin{abstract}                % Abstract of not more than 250 words.
The paper proposes a framework for the construction of solutions to a hybrid dynamical system that exhibit Zeno behavior. A new approach that enables solution to be prolonged after reaching its Zeno time is developed. It allows for a comprehensive stability analysis and asymptotic behavior characterization of such solutions. The results are applicable to a wide class of hybrid systems and match with practical experience of simulation of real-world phenomena. Moreover they are potentially useful for applications to interconnections of hybrid systems.
\end{abstract}


\end{frontmatter}
%===============================================================================

\section{Introduction}\label{ss1}

Processes that combine continuous and discontinuous behavior naturally arise in a variety of real-world applications such as robotics, biological systems, chemical kinetics, logistics and networked control systems. The basic framework to model and analyse such a behavior is impulsive differential equations~\cite{samoilenko1987differential,lakshmikantham1989theory,samoilenko1995impulsive}.
%Theory of impulsive differential equations dates back to the papers of Milman and Myshkis~\cite{milman1960stability}, Myshkis and Samoilenko~\cite{myshkis1967systems} that appeared in 1960-th. In these works a concept of solution to impulsive differential equation, conditions for its existence and uniqueness, fundamental stability properties have been introduced. Later in 1970-th a rigorous theory of impulsive differential equations with fixed and non-fixed moments of impulsive actions has been developed by Samoilenko, Perestyuk and their descendants. Their results include classification of impulsive differential equations depending on the character of impulsive actions~\cite{samoilenko1987differential,samoilenko1995impulsive}, stability characterization of solutions~\cite{Sam1977,samoilenko1982stability}, existence of periodic and almost-periodic solutions~\cite{samoilenko1978periodic,PSh79,akhmetov1983almost} and invariant sets~\cite{Per84}, Lyapunov-like theorems and associated dwell-time conditions~\cite{gurgula1982second,gurgula1982stability}, justification of averaging method~\cite{samoilenko1961application,Sam1967,mitropol1985averaging}.
Besides this theory we would like to mention the other more recent developments in the related fields like hybrid dynamical systems~\cite{van2000introduction}, dynamic equations on time scales~\cite{bohner2012dynamic}, discontinuous dynamical systems~\cite{akhmet2010principles}, switched systems~\cite{liberzon2012switching} and hybrid automata~\cite{henzinger2000theory}.

Throughout the paper we will use one of the most recent and rapidly developing framework --- a hybrid dynamical system proposed in~\cite{goebel2012hybrid}.
This framework is one of the most general and includes a majority of other classes of systems that model processes with continuous and discontinuous behavior. Moreover a variety of novel results are developed in \cite{goebel2012hybrid} that are not available in the other frameworks. Also this framework appears well-adapted to the control-related problems. In particular the introduction of input-to-state stability (ISS) concept for hybrid dynamical systems gave a strong push and motivated a fast development of new methods for stability analysis of hybrid systems with exogenous input~\cite{CaT09}. The questions on robustness of ISS for hybrid systems were considered in~\cite{CaT13}. In recent years a considerable attention is paid to the stability analysis of interconnections of hybrid dynamical systems. Small-gain approach proved to be an effective tool for stability analysis of solutions to interconnections and networks of a large scale~\cite{sanfelice2011interconnections,dashkovskiy2013input,MYL14,San14,LNT14}. In spite of these developments, interconnections of hybrid systems are considered only under strong constraints, that are often not compatible with applications~\cite{sanfelice2011interconnections,dashkovskiy2013hybrid}.

The simplest example of unsolved problem is related to a bouncing ball modelled by hybrid dynamical system. The origin of the bouncing ball system is in some sense asymptotically stable (for a precise definition see Definition~\ref{dugpas}). There is a variety of Lyapunov-like theorems in \cite{goebel2012hybrid} to verify this. However if we consider two such balls as one system (a so-called vacuous interconnection) then there are no methods to prove the asymptotic stability of the origin for the entire system. Moreover, this system is not asymptotically stable in the framework of \cite{goebel2012hybrid}. It seams that this framework is not suitable for modeling this rather simple mechanical system, however we claim that the stability problem can be resolved by a minor extension of the theory developed in \cite{goebel2012hybrid}. For more details of the just mentioned problem we refer to Section~\ref{motivation} and Figure~\ref{mynewimage}. Here we only mention that this problem is caused by the Zeno solutions characterized by infinitely many impulsive jumps over a finite period of time. Such solutions are not defined after this time period. %This leads to a problem when this time period is different for each of the two bouncing balls.

Several approaches were proposed to cope with this problem. Some of these methods enable a solution to be prolonged beyond its Zeno time but only for certain classes of hybrid systems. In~\cite{johansson1999regularization}, a so-called regularization technique has been proposed and was illustrated for particular examples. It is based on perturbing the hybrid system in order to obtain non-Zeno solution, and then taking the limit as the perturbation goes to zero. A more formal procedure for obtaining generalized solutions of Zeno hybrid system via regularization was presented in \cite{goebel2004hybrid,sanfelice2008generalized}.
%Following similar guidelines, generalized solutions for unilaterally constrained mechanical systems are analyzed in \cite{miller2007generalized}.
For a particular class of Lagrangian hybrid systems, a solution switches to a holonomically constrained dynamical system after the Zeno point is reached \cite{ames2006there}, \cite{or2011stability}. In the closely related class of switched systems \cite{liberzon2012switching}, \cite{shorten2007stability}, a solution may converge to a switching surface in a finite time, along with increasingly fast switching events near this surface. This phenomenon is called chattering. In this case, the solutions can be extended by considering the set-valued Filippov solution \cite{filippov1988differential}, which involves sliding along the switching surface. In \cite{cuijpers2001beyond}, a solution prolongation beyond Zeno was proposed by introducing the concept of transition over infinite sequence and accumulation-closed transition systems. Finally, considerable achievements were made from a computer science viewpoint. The existence of infinitely many discrete events over a finite period of time force simulators to ignore some events or looping indefinitely. The ways to overcome these problems were proposed in \cite{Kon2016} by introducing new algorithms for event detection and localization.

A peculiarity of hybrid dynamical systems is that the concept of time is characterized by two parameters: the amount of time passed and the number of jumps that have occurred. According to~\cite{goebel2012hybrid}, a certain subset of $\R_{\geq 0}\times\N_0$ is called hybrid time domain. More general rules for constructing hybrid time domains were proposed in \cite{collins06,davoren2008topologies}. In \cite{collins06}, the concept of generalized hybrid time domain has been introduced where a discrete-time axis was generalised to a countable ordinal that can have infinitely many accumulation points, which correspond to Zeno occurrences. This approach enables to prolong solutions to a hybrid system beyond Zeno time. However in \cite{collins06,davoren2008topologies} authors did not study stability properties of prolonged solutions which is in the main focus of our paper.

The aim of this paper is to develop an approach for solutions prolongation over the Zeno time and to study their stability properties. In view of the well developed stability theory in \cite{goebel2012hybrid} we aim to introduce some minor extensions in its framework so that we still can use results from \cite{goebel2012hybrid}. For this reason we do not follow such deep modifications of the hybrid time domain notion as in \cite{collins06,davoren2008topologies}. Our slight extension enables stability analysis of hybrid systems beyond Zeno points.% which is the second aim of this paper. As the third contribution we would like to mention a rather wide and up to date literature review on the Zeno solutions to hybrid systems, their prolongation and related stability problems.

In this paper we propose an approach to extend a solution to a hybrid dynamical system beyond its Zeno time without destroying the key concepts of \cite{goebel2012hybrid}. In our mind, a natural way is to prolong Zeno solution from its $\omega$-limit point. For this purpose we adapt hybrid framework from \cite{goebel2012hybrid} by introducing a three dimensional hybrid time domain and redefining the concept of solution.% We believe that our approach is applicable to a wide class of processes that exhibit Zeno behavior and illustrate it with the examples.

The rest of the paper is organized as follows. In Section 2 we recall some basic definitions from the theory of hybrid dynamical systems. A motivating example is given in Section 3. A new approach for solution construction is presented in Section 4. In Section 5 we prove a series of propositions that enable stability analysis of solutions to a hybrid dynamical system with Zeno behavior. An illustrative example is given there. A short discussion on open problems in Section 6 completes the paper.

\section{Preliminary notion and definitions}

The following notation and definitions are taken from \cite{goebel2012hybrid}:
\begin{equation}\tag{$\H$}\label{Hybr}
\begin{cases}
\dot x &= \quad f(x),\quad x \in C, \\
x^+ &= \quad g(x), \quad x \in D.
\end{cases}
\end{equation}
The state $x\in\R^n$, $n\in\N$ can change according to the differential equation $\dot x = f(x)$ while $x \in C$, and it can change according to the difference equation $x^+ = g(x)$ while $x \in D$. The sets $C\subset\R^{n}$ and $D\subset\R^{n}$ are called the flow and the jumps sets respectively, functions $f:C\to\R^{n}$ and $g:D\to\R^{n}$ are the flow and jump maps. The data of the hybrid system \ref{Hybr} is given by $(C,f,D,g)$.

The parametrization of a solution to the hybrid system \ref{Hybr} is given by two parameters: $t\in\R_{\geq 0}=[0,\infty)$, the amount of time passed, and $j\in\N_0=\N\cup\{0\}$, the number of jumps that have occurred. A certain subset of $\R_{\geq 0}\times\N_0$ can correspond to evolutions of hybrid systems. Such sets are called hybrid time domains.
\begin{definition}[Hybrid time domain]\label{htd}
Let $t_0\leq t_1 \leq t_2 \leq t_3 \leq \ldots$. A subset
\begin{equation*}
E=\bigcup_j([t_j,t_{j+1}],j) \subset\R_{\geq 0}\times \N_0
\end{equation*}
is a hybrid time domain if it is a union of a finite or infinite sequence of intervals $[t_j,t_{j+1}]\times\{j\}$, with the last interval (if existent) possibly of the form $[t_j,T)$ with $T$ finite or $T=\infty$.
\end{definition}

Given a hybrid time domain $E$ we denote:
\begin{equation*}
\begin{split}
&\sup_t{E}=\sup\{t\in\R_{\geq 0}: \exists j\in\N_0\text{~such~that~}(t,j)\in E\}, \\
&\sup_j{E}=\sup\{j\in\N_{0}: \exists t\in\R_{\geq 0}\text{~such~that~}(t,j)\in E\}. %, \\
\end{split}
\end{equation*}

\begin{definition}[Hybrid arc]\label{darc}
A function $\phi:E\to\R^n$ is a hybrid arc if $E$ is a hybrid time domain and if for each $j\in\N$, the function $t\to\phi(t,j)$ is locally absolutely continuous on the interval $I^j=\{t:(t,j)\in E\}$.
\end{definition}
Given a hybrid arc $\phi$, the notation $\dom\phi$ represents its domain, which is a hybrid time domain.

\begin{definition}[Complete hybrid arc]
A hybrid arc $\phi:E\to\R^n$ is called complete if $\dom\phi$ is unbounded, i.e., if $\sup_t{E} + \sup_j{E}=\infty$.
\end{definition}

\begin{definition}[Zeno hybrid arc]
A hybrid arc $\phi:E\to\R^n$ is called Zeno if it is complete and $\sup_{t}\dom\phi<\infty$.
\end{definition}
The existence of a Zeno hybrid arc means that an infinite number of jumps occurs during a finite time. The time $\tau=\sup_{t}{\dom\phi}$ is called a Zeno time.

\begin{definition}[Solution to a hybrid system]\label{dsol}
A hybrid arc $\phi$ is a solution to the hybrid system $\H$ if $\phi(0,0)\in\bar C \cup D$ and
\begin{itemize}
\item[(S1)] for all $j\in\N$ such that $I^j:=\{t:(t,j)\in \dom\phi\}$ has nonempty interior
\begin{equation*}
\begin{array}{lll}
\phi(t,j)\in C &\text{~for all} &t\in \interior I^j,\\
\dot\phi(t,j)=f(\phi(t,j)) &\text{for almost all} &t \in I^j;
\end{array}
\end{equation*}
\item[(S2)] for all $(t,j)\in\dom\phi$ such that $(t,j+1)\in\dom\phi$,
\begin{equation*}
\begin{array}{ll}
\phi(t,j)\in D, & \phi(t,j+1)=g(\phi(t,j)).
\end{array}
\end{equation*}
\end{itemize}
\end{definition}

The properties of hybrid arcs (like completeness, Zeno, etc.) are automatically extended on the corresponding solutions.

\begin{definition}[Maximal solution]
A solution $\phi$ to $\H$ is maximal if there does not exist another solution $\psi$ to $\H$ such that $\dom\phi$ is a proper subset of $\dom\psi$ and $\phi(t,j)=\psi(t,j)$ for all $(t,j)\in\dom\phi$.
\end{definition}
Let $\mathcal S_{\H}(\A)$ denote the set of all maximal solutions $\phi$ to a hybrid system $\H$ with $\phi(0,0)\in\A$.

\begin{definition}[Strong forward pre-invariance]
A set $\A\subset\R^n$ is said to be strongly forward pre-invariant (SFpI) if for every $\phi\in\mathcal S_{\H}(\A)$, $\rge\phi\subset \A$, where $\rge\phi=\{y\in\R^n: \exists (t,j)\in\dom\phi \text{~such~that~} y=\phi(t,j)\}$.
\end{definition}

For a precise definition of stability we recall the definitions of standard functions and distance to a closed set. A function $\alpha:\R_{\geq0}\to \R_{\geq0}$ is called a class-$\K_\infty$ function ($\alpha\in\K_\infty$) if $\alpha$ is zero at zero, continuous, strictly increasing, and unbounded. A function $\rho:\R_{\geq 0}\to \R_{\geq 0}$ is positive definite ($\rho\in\PD$) if $\rho(s)>0$ for all $s>0$ and $\rho(0)=0$. Given a vector $x\in\R^n$ and a closed set $\A\subset\R^n$, the distance of $x$ to $\A$ is defined by $|x|_{\A}:=\inf\limits_{y\in\A}|x-y|$.

\begin{definition}[Uniform global pre-asymptotic stability]\label{dugpas}
Let $\A\subset\R^n$ be closed. The set $\A$ is said to be
\begin{itemize}
\item uniformly globally stable (UGS) if there exists a function $\alpha\in\K_\infty$ such that any solution $\phi$ to $\H$ satisfies $|\phi(t,j)|_{\A}\leq \alpha(|\phi(0,0)|_{\A})$ for all $(t,j)\in\dom\phi$;
\item uniformly globally pre-attractive (UGpA) if for each $\varepsilon>0$ and $r>0$ there exists $T>0$ such that, for any solution $\phi$ to $\H$ with $|\phi(0,0)|_{\A}\leq r$, $(t,j)\in\dom\phi$ and $t+j \geq T$ imply $|\phi(t,j)|_{\A}\leq \varepsilon$;
\item uniformly globally pre-asymptotically stable (UGpAS) if it is both uniformly globally stable and uniformly globally attractive.
\end{itemize}
\end{definition}

\begin{definition}[$\omega$-limit set of a hybrid arc]
The $\omega$-limit set of a hybrid arc $\phi:\dom\phi\to\R^n$, denoted $\Omega(\phi)$, is the set of all points $x\in\R^n$ for which there exists a sequence $\{(t,j)_i\}^{\infty}_{i=1}$ of points $(t_i,j_i)\in\dom\phi$ with $\lim\limits_{i\to\infty}{t_i+j_i}=\infty$ and $\lim\limits_{i\to\infty}\phi(t_i,j_i)=x$. Every such point is an $\omega$-limit point of $\phi$.
\end{definition}

\section{Motivating example}\label{motivation}
Consider two hybrid dynamical systems \ref{H} with states $x_i\in\R^{n_i}$ and inputs $u_i\in\mathcal U_i\subset\R^{m_i}$
\begin{equation}\tag{$\H_{i}$}\label{H}
\begin{cases}
\dot x_{i} &= \quad f_{i}(x_{i},u_{i}),\quad (x_{i},u_{i})\in C_{i}, \\
x_{i}^+ &= \quad g_{i}(x_{i},u_{i}), \quad (x_{i},u_{i}) \in D_{i},
\end{cases}
\end{equation}
where $n_i, m_i\in\N$, $i=1,2$. The sets $C_i\subset\R^{n_i}\times {~\mathcal U}_i$ and $D_i\subset\R^{n_i}\times {~\mathcal U}_i$ define the flow and the jumps sets respectively, functions $f_i:C_i\to\R^{n_i}$ and $g_i:D_i\to\R^{n_i}$ are the flow and jump maps. The data of the hybrid system \ref{H} is given by $(C_i,f_i,D_i,g_i)$.

Let us interconnect these two systems with $u_1=h_1(x_1)$ and $u_2=h_2(x_2)$, where functions $h_1:\R^{n_1}\to\mathcal U_2$, $h_2:\R^{n_2}\to\mathcal U_1$. Then the entire interconnection can be represented as a single hybrid dynamical system $\H$ with data $(C,f,D,g)$, where its state is $x:=(x_1,x_2)\in\R^{n_1}\times\R^{n_2}$, its flow set is $$C:=\{x:(x_1,h_2(x_2))\in C_1\}\cap\{x:(x_2,h_1(x_1))\in C_2\},$$ its flow map is $f(x):=(f_1(x_1,h_2(x_2)),f_2(x_2,h_1(x_1)))$, its jump set is $$D:=\{x:(x_1,h_2(x_2))\in D_1\}\cup\{x:(x_2,h_1(x_1))\in D_2\}$$ and its jump map is $g(x):=(\tilde g_1(x_1,h_2(x_2)),\tilde g_2(x_2,h_1(x_1)))$ with
\begin{equation*}
\tilde g_1(x):=\begin{cases}
g_1(x_1,h_2(x_2)),\quad\text{if}\quad (x_1,h_2(x_2))\in D_1,\\
x_1\quad\quad\quad\quad\quad\quad\text{otherwise,}
\end{cases}
\end{equation*}
\begin{equation*}
\tilde g_2(x):=\begin{cases}
g_2(x_2,h_1(x_1)),\quad\text{if}\quad (x_2,h_1(x_1))\in D_2,\\
x_2\quad\quad\quad\quad\quad\quad\text{otherwise.}
\end{cases}
\end{equation*}

\begin{figure}[h]
  \center
  \includegraphics[width=0.5\textwidth]{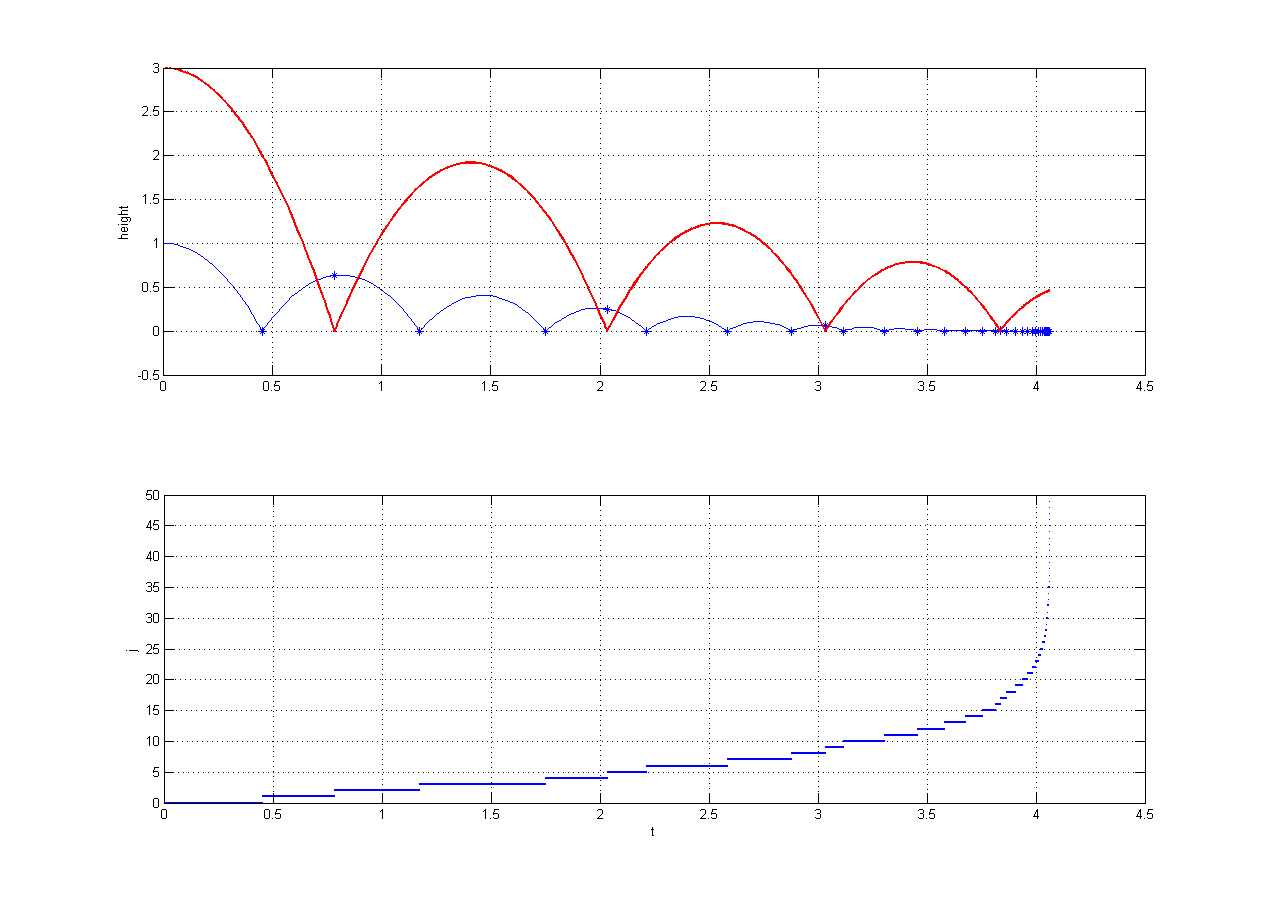}
  \caption{Evolution of the height coordinates (top) and hybrid time domain (bottom) of two vacuously interconnected bouncing balls started with the height 3 (red) and 1 (blue) respectively. The state of the system converges to some point away from the origin when hybrid time goes to infinity ($t+j\to\infty$). The solution is not defined beyond the observed Zeno time.}
\label{mynewimage}
\end{figure}

In the literature~\cite{dashkovskiy2013input,dashkovskiy2013hybrid} such choice of the flow set $C$ and the jump set $D$ is called natural. An important fact is that an interconnection of two hybrid systems $\H_{1}$ and $\H_{2}$ is a hybrid system of the form $\H$. So one may use a variety of previously developed methods and techniques (for instance from \cite{goebel2012hybrid}) for a qualitative characterization of solutions and the problem of a comprehensive analysis of interconnections seems to be solved. However an essential problem appears in this context. It was discussed in \cite{sanfelice2011interconnections} and caused by the interconnection of a hybrid system with Zeno solution and a hybrid system with continuous complete solution. Such interconnection has a Zeno solution that is not a part of the set of solutions to every subsystem. Another good illustration of this problem is a vacuous interconnection of several bouncing balls when the balls start from different initial positions. The solution of such model may not allow all the balls to reach their own Zeno time as the original model of each bouncing ball does (see Figure~\ref{mynewimage}). This leads to unnatural loss of asymptotic stability of the origin.

In this paper we propose a way to extend the hybrid framework~\cite{goebel2012hybrid} in order to cope with aforementioned problems.

\section{Hybrid framework extension}
The main source of the problems stated in the motivation section is that a solution to a hybrid system is not defined beyond its Zeno time. However some experiments from real life like bouncing ball argue that a solution should be prolonged over its Zeno time. A bouncing ball after reaching the resting state continues to lie while time is counting further and further. This motivates us to allow solution to continue its evolution after reaching Zeno time. In our extended framework, Zeno solution continues its evolution from an $\omega$-limit point after reaching its Zeno time. It enables us to construct solutions that reflect real-world observations and to perform their stability analysis.

To describe the evolution of solution to a hybrid system we introduce a new notion of hybrid time domain. It tracks not only the elapsed time and the number of impulsive jumps, but also the number of Zeno points occurred during the evolution process. Similar to a classical hybrid time domain from Definition~\ref{htd}, only certain subsets of $\R_{\geq 0}\times \N_0 \times \N_0$ can correspond to evolutions of hybrid systems.

\begin{definition}[Extended hybrid time domain]
Let $\{t_{j,k}\}$ be a set of time moments such that $t_{j,k}\leq t_{j+1,k}$ $\forall j,k\in\N_0$ and $t_{j,k}\leq t_{i,k+1}$ $\forall i,j,k\in\N_0$. A subset
\begin{equation*}
\tilde E=\bigcup_{j,k}([t_{j,k},t_{j+1,k}],j,k) \subset\R_{\geq 0}\times\N_0\times\N_0
\end{equation*}
is an extended hybrid time domain if it is a union of a finite or infinite set of intervals $[t_{j,k},t_{j+1,k}]\times\{j\}\times\{k\}$, with the last interval (if existent) possibly of the form $[t_{j,k},T)$ with $T$ finite or $T=\infty$.
\end{definition}
Index $k$ corresponds to the number of encountered Zeno behaviors. For a given extended hybrid time domain $\tilde E$ we denote:
\begin{equation*}
\begin{split}
\sup_{Zeno}{\tilde E}=\sup\{k\in\N_{0}: \exists t\in\R_{\geq 0}, j\in\N_0\text{~s.t.~}(t,j,k)\in \tilde E\}. %, \\
\end{split}
\end{equation*}
Note that for any extended hybrid time domain we can fix an admissible index $k$ and consider ist subset corresponding to this $k$. Its projection onto
$\R_{\geq 0}\times\N_0$ (defined by dropping $k$) is the "classical" hybrid time domain from Definition~\ref{htd}.

An \emph{extended solution} $\tilde\phi$ is a function defined on an \emph{extended hybrid time domain}. Before reaching the first Zeno time the extended solution $\tilde\phi$ coincides with the ''classical'' solution $\phi$ to a hybrid system: $\tilde\phi(t,j,0)\equiv \phi(t,j)$ for all $(t,j)\in\dom\phi$.

In the original framework~\cite{goebel2012hybrid} a state $x(t)=\varphi(t,j)\in\R^n$ can evolve along a trajectory of differential equation $\dot x = f(x)$ while $x(t)\in C$. At the time $s\in\R_{\geq 0}$ when $x(s)\in D$ it can be instantly transferred into a new position $g(x(s))$ and the value of the corresponding jump index in hybrid time domain increases by 1 so that $\varphi(s,j+1)=g(x(s))$. In our settings we add one more rule to construct solution to a hybrid dynamical system:
\begin{itemize}
\item if for some fixed $k\in\N_0$ hybrid arc $\tilde\phi(t,j,k)$ is Zeno with non-empty $\omega$-limit set then a solution $\tilde\phi$ to a hybrid system $\H$ is prolonged with initial condition $\tilde\phi(\tau,0,k+1)\in\Omega(\phi)$, where $\tau$ is the Zeno time for hybrid arc $\tilde\phi(\cdot,\cdot,k)$.
\end{itemize}

Our extended solution $\tilde\phi$ is a concatenation of classical hybrid arcs $\phi^i(t,j)$ with initial conditions $\phi^0(0,0)=\xi$, $\phi^1(\tau_1,0)\in\Omega(\phi^0)$, $\phi^2(\tau_2,0)\in\Omega(\phi^1)$ and so on, where $\tau_i$ is the Zeno time for the hybrid arc $\phi^{i-1}$.

A new rule of extended solution's construction leads to the following properties of the corresponding extended hybrid time domain $\tilde E$ : if the point $(t_{0,k+1},0,k+1)\in\tilde E$ then there exist infinitely many points of the form $(t_{\cdot,k},\cdot,k)\in\tilde E$ such that $\lim\limits_{j\to\infty}t_{j,k}=t_{0,k+1}$.

In general, an $\omega$-limit set $\Omega(\phi)$ may consist of several or infinitely many points. According to our new rule a single initial point can generate multiple solutions. Such situation appears, for example, in modelling of water tanks system (see \cite{alur1997modularity} for details). This system has Zeno arcs with two $\omega$-limit points. Therefore two different extended solutions will be generated from a single initial point. This is quite natural since the considered physical process can evolve according to both of solutions in a real experiment.

\begin{remark}
For a given Zeno hybrid arc $\varphi$ the $\omega$-limit point $x\in\Omega(\varphi)$ is a limit of a sequence of points $\{x_i\}$ of the state space such that $x_i\in D$, $i=1,2,\ldots$. A behavior of the corresponding extended solutions now heavily depends on the properties of the jump set $D$. If $D$ is a closed set (it means that it contains all its limit points) then the extended solution continue its evolution from limit point $x\in D$ and therefore should jump. This can lead to eventually discrete solution. If the jump map $D$ is an open set then the extended solution can continue its evolution from the limit point $x\in C$ and therefore can be prolonged beyond Zeno time of the ordinary time axis. From this viewpoint it is quite natural to model real processes with an open jump set $D$ as it is done in the following example.
\end{remark}

%The following example demonstrates the proposed approach.

\begin{example}\label{ex1}
Consider a vacuous interconnection of two bouncing balls. Let $x_1, x_3 \in \R_{\geq 0}$ stand for the heights of the balls and $x_2, x_4 \in \R$ stand for the corresponding velocities. Then system has the form
\begin{equation*}
\begin{array}{lcl}
\dot x_1 & = & x_2, \\
\dot x_2 & = & -\gamma(x_1,x_2),
\end{array}
\quad
\begin{array}{lcl}
\dot x_3 & = & x_4, \\
\dot x_4 & = & -\gamma(x_3,x_4),
\end{array}
\quad x\in C,
\end{equation*}

\begin{equation*}
\begin{gathered}
\begin{array}{lcl}
{x_1}^+ & = & x_1, \\
{x_2}^+ & = & \begin{cases}-\lambda x_2, x\in D_1,\\ x_2,\quad\text{~} x\not\in D_1, \end{cases}
\end{array}\hspace{-3mm}
\begin{array}{lcl}
{x_3}^+ & = & x_3, \\
{x_4}^+ & = & \begin{cases}-\lambda x_4, x\in D_2,\\ x_4,\quad\text{~} x\not\in D_2, \end{cases}
\end{array}
\hspace{-2mm}x\in D,
\end{gathered}
\end{equation*}

\begin{equation*}
\begin{split}
C_1&=\{x\in\R^4: x_1>0 \text{~or~} x_1=0, x_2\geq0\}, \\
D_1&=\{x\in\R^4: x_1=0, x_2<0\},\\
C_2&=\{x\in\R^4: x_3>0 \text{~or~} x_3=0, x_4\geq0\}, \\
D_2&=\{x\in\R^4: x_3=0, x_4<0\},\\
C&=C_1\cap C_2, \quad D=D_1\cup D_2,
\end{split}
\end{equation*}

where $\lambda\in(0,1)$ is the restitution coefficient, $\gamma:\R^2\to\R$ is given by
\begin{equation}\label{eqgamma}
\gamma(a,b)=
\begin{cases}
0,\quad \text{~if~} a=b=0,\\
9,81\quad\text{~otherwise}.
\end{cases}
\end{equation}
\begin{figure}[h]
  \includegraphics[width=0.5\textwidth]{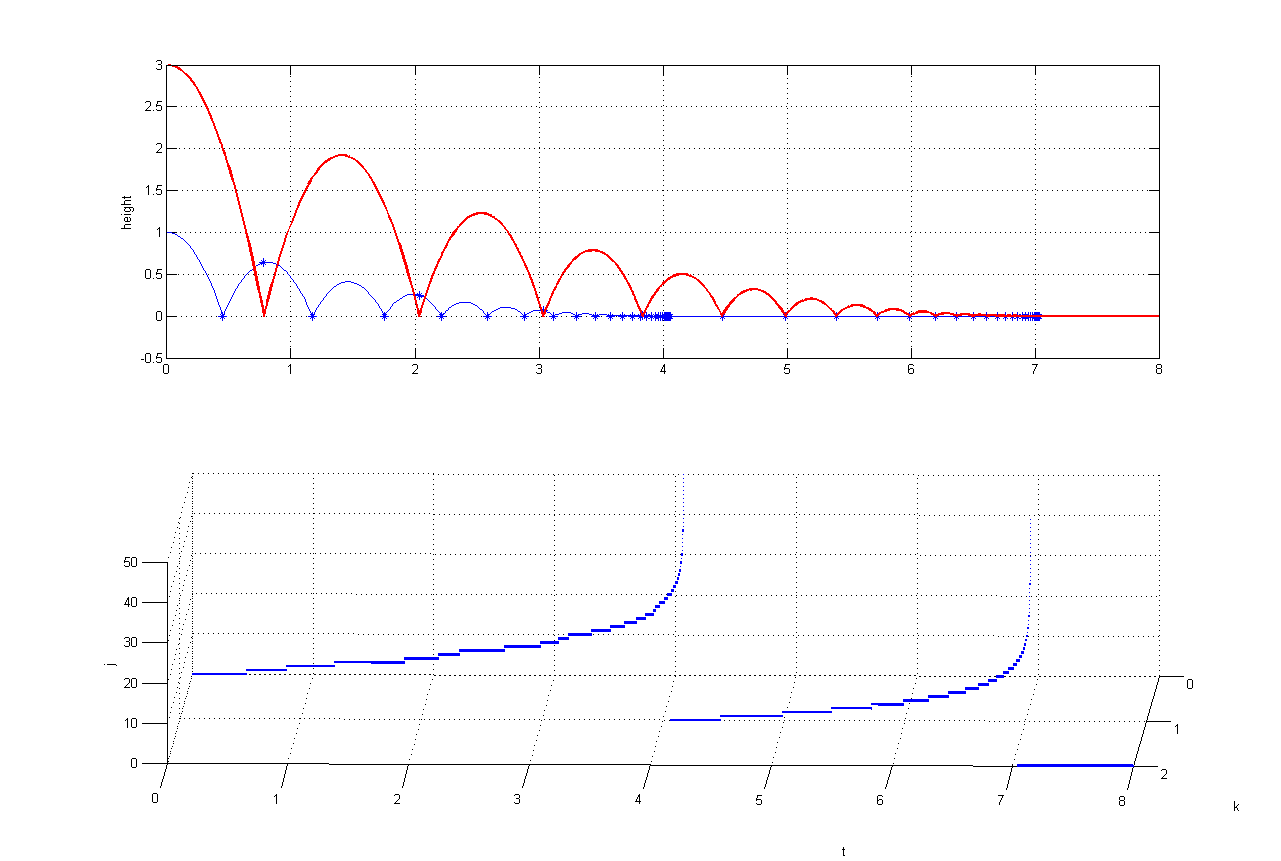}
  \caption{Evolution of the height coordinates (top) and extended hybrid time domain (bottom) of two vacuously interconnected bouncing balls started with the height 3 (red) and 1 (blue) respectively.}
\label{image2}
\end{figure}
\end{example}
A numerical simulation is presented on Figure~\ref{image2}. The arc that corresponds to Zeno index $k=0$ fully coincides with the one from the original framework \cite{goebel2012hybrid}.
Its $\omega$-limit set consists of a single point (in this case the uniqueness of solution is preserved). At the Zeno time of the blue ball the solution is now prolonged from this point and Zeno index is increased by $1$.
The extended solution exhibits a further Zeno behavior and its $\omega$-limit set is just the origin.
At the Zeno time of the red ball, the solution is prolonged from its $\omega$-limit set $(0,0,0,0)\in\R^4$ which is a single point again. The last arc of solution is trivial and purely continuous with $\sup_{t}\dom\phi=\infty$. The concatenated solution corresponds to our experience.

\section{Stability analysis}\label{sstab}
In this section we introduce a new stability notion in order to describe asymptotic behavior of extended solutions to hybrid systems. Two auxiliary lemmas will be needed to justify stability characterization.

\begin{lemma}\label{sfpi}
UGS of a set $\A$ implies its strong forward pre-invariance.
\end{lemma}
\begin{pf}
Suppose it is not true. Let there exist a solution $\phi$ to $\H$ with $\phi(0,0)\in\A$ and a point $x^*\in\rge\phi(t,j)$ such that $x^*\not\in\A$ for some $(t,j)\in\dom\phi$. It means that there exists $\delta>0$ such that $|\phi(t,j)|_{\A}=\delta$. Then from Definition~\ref{dugpas} it follows directly that there is a function $\alpha\in\K_{\infty}$ such that
\begin{equation*}
0<\delta=|\phi(t,j)|_{\A}\leq\alpha(|\phi(0,0)|_{\A})=\alpha(0)=0.
\end{equation*}
The contradiction proves that every solution $\phi$ starting in $\A$ remains in this set: $\phi(t,j)\in\A$ for all $(t,j)\in\dom\phi$.
\end{pf}

\begin{lemma}\label{olimit}
Let $\A$ be UGpAS and every arc $\phi$ with initial condition in $\{C\cup D\} \setminus \A$ have a non-empty $\omega$-limit set, then $\Omega(\phi) \subset \A$.
\end{lemma}
\begin{pf}
Suppose it is not true. Consider a solution $\phi^*$ to $\H$ with $|\phi^*(0,0)|_{\A}\leq r$ that has an $\omega$-limit point $\xi$ outside the set $\A$. It means that there exists a sequence $\{(t,j)_i\}^{\infty}_{i=1}$ of points $(t_i,j_i)\in\dom\phi$ with $\lim\limits_{i\to\infty}{t_i+j_i}=\infty$ and
\begin{equation}\label{lp}
\lim\limits_{i\to\infty}\phi^*(t_i,j_i)=\xi\not\in\A.
\end{equation}
Then there exists $\delta>0$ such that $|\xi|_{\A}=\delta$. From the UGpAS of the set $\A$ it follows that for $\varepsilon=\frac{\delta}{2}$ there exists $T>0$ such that for every $(t,j)\in\dom\phi$ with $t+j\geq T$ follows $|\phi(t,j)|_{\A}\leq \frac{\delta}{2}$.

However the existence of the limit (\ref{lp}) guarantees that for any $d>0$ there exist a $d$-neighbourhood $U_d(\xi)$ and $(t^*, j^*)\in\dom\phi$ such that $\phi(t^*,j^*)\in U_d(\xi)$. Choosing $d$ small enough to satisfy the conditions $t^* + j^* \geq T$ and $U_d(\xi)\cap U_{\frac{\delta}{2}}(\A)=\emptyset$ leads to $|\phi(t^*,j^*)|_{\A} > \frac{\delta}{2}$, which contradicts the UGpAS of the set $\A$. This proves that the $\omega$-limit point $\xi\in\A$.
\end{pf}

\begin{definition}[$\H\cap\A$]\label{phase}
If $\A\subset C\cup D$ is UGpAS for $\H$, then $\A$ can be considered as the state space for a new hybrid system with the new flow set $C\cap\A$ and the new jump set $D\cap\A$. We will denote this new system by $\H\cap\A$.
\end{definition}
Indeed, from Lemma \ref{sfpi}, UGpAS implies SFpI of the set $\A$ so every solution with initial condition in $\A$ will remain there for all $(t,j)\in\dom\phi$. In the case when $\phi$ is a Zeno solution it will be prolonged from a point of its $\omega$-limit set $\Omega(\phi)$. From Lemma \ref{olimit} follows that $\Omega(\phi)\subset \A$, so the solution will again remain in the set $\A$. It means that extended solutions of the system $\H$ with initial conditions in $\A$ will coincide with extended solutions of the system $\H\cap\A$ with the corresponding initial conditions. Since $\A\subset C\cup D$, no new solution will be generated.

For a comprehensive description of asymptotic behavior of extended solutions we introduce a new definition of stability over Zeno.

\begin{definition}[UGpASoZ]\label{UGpASoZ}
Let $\A\in\R^n$ be closed. The set $\A$ is said to be
\begin{itemize}
\item[(i)] uniformly globally stable over Zeno (UGSoZ) if there exists a function $\alpha\in\mathcal K_\infty$ such that any solution $\tilde\phi$ to $\H$ satisfies $|\tilde\phi(t,j,k)|_{\mathcal A}\leq \alpha(|\tilde\phi(0,0,0)|_{\mathcal A})$ for all $(t,j,k)\in\dom\tilde\phi$;
\item[(ii)] uniformly globally pre-attractive over Zeno (UGpAoZ) if for each $\varepsilon>0$ and $r>0$ there exist $T>0$ and $K\geq 0$ such that, for any solution $\tilde\phi$ to $\H$ with $|\tilde\phi(0,0,0)|_{\mathcal A}\leq r$, from $(t,j,k)\in\dom\tilde\phi$ with either $t+j \geq T$, $k=K$ or $k>K$ or $t+j \geq T$, $k=\sup_{Zeno}{\dom\tilde\varphi}$, $\sup_{Zeno}{\dom\tilde\varphi}<K$ it follows that $|\tilde\phi(t,j,k)|_{\mathcal A}\leq \varepsilon$;
\item[(iii)] globally pre-attractive over Zeno (GpAoZ) if for each $\varepsilon>0$, $r>0$, and for any solution $\tilde\phi$ to $\H$ with $|\tilde\phi(0,0,0)|_{\mathcal A}\leq r$, there exist $T>0$ and $K\geq 0$ such that from $(t,j,k)\in\dom\tilde\phi$ with either $t+j \geq T$, $k=K$ or $k>K$ it follows that $|\tilde\phi(t,j,k)|_{\mathcal A}\leq \varepsilon$;
\item[(iv)] uniformly globally pre-asymptotically stable over Zeno (UGpASoZ) if it is both UGSoZ and UGpAoZ.
\end{itemize}
\end{definition}

The conditions for the pre-attractivity actually mean that all solutions will reach the $\varepsilon$-neighbourhood of the set $\A$ no later than at the time $T$ after $K$-th Zeno occurrence. The uniformity means that $T$ and $K$ are the same for all solutions. If a solution does not undergo such number ($K$) of Zeno occurrences then it should reach the corresponding $\varepsilon$-neighbourhood no later than time $T$ after its last Zeno.

\begin{thm}\label{main}
Let there exist a finite sequence $\A_n\subset \A_{n-1}\subset\ldots\subset \A_1\subset \A_0= C\cup D$ such that $\A_i$ is UGpAS for the system $\H\cap\A_{i-1}$, $i=1,\ldots,n$ and for all initial values $\phi(0,0)\in \A_{i-1} \setminus \A_i$, $i=1,\ldots,n-1$ the corresponding solutions $\phi$ are Zeno with non-empty $\omega$-limit sets. Then $\A_n$ is UGpASoZ.
\end{thm}
\begin{pf}
If $n=1$ then UGpAS of the set $\A_1$ implies its UGpASoZ with $K=0$. Let us consider the case $n=2$. First we prove stability of the set $\A_2$. From UGS of the sets $\A_1$ and $\A_2$ it follows that there exist functions $\alpha_1,\alpha_2\in \K_\infty$ such that
\begin{equation*}
\begin{split}
|\varphi(t,j)|_{\A_1}&\leq \alpha_1(|\phi(0,0)|_{\mathcal A_1})\quad \forall\varphi(0,0)\in A_0, \\
|\varphi(t,j)|_{\A_2}&\leq \alpha_2(|\phi(0,0)|_{\mathcal A_2})\quad \forall\varphi(0,0)\in A_1
\end{split}
\end{equation*}
and for all $(t,j)\in\dom\phi$. Then the extended solution $\tilde \varphi$ satisfies
\begin{equation*}
\begin{split}
|\tilde\varphi(t,j,k)|_{\A_2}&\leq \alpha_{12}(|\tilde\phi(0,0,0)|_{\mathcal A_2})\quad \forall\varphi(0,0)\in A_0
\end{split}
\end{equation*}
for all $(t,j,k)\in\dom\tilde\phi$ with $\alpha_{12}(s)=\max\{\alpha_{1}(s), \alpha_{2}(s)\}$. Hence $\A_2$ is UGSoZ.

Now we will prove the pre-attractivity. For this purpose we will show that each extended solution $\tilde\phi$ to the hybrid system $\H$ satisfies the uniform pre-attractivity conditions (\emph{ii}) of Definition~\ref{UGpASoZ} with respect to the set $\A_2$. Note that since $\A_2$ is UGpAS for the system $\H\cap\A_1$ every solution initiated from $\A_1$ satisfies the pre-attractivity conditions of the Definition \ref{dugpas}: $\forall \varepsilon,r>0$ there exists $T_1>0$ such that any solution $\varphi$ to $\H\cap\A_1$ such that $|\phi(0,0)|_{\A_2}\leq r$, $\phi(0,0)\in\A_1$ satisfies $|\phi(t,j)|_{\A_2}\leq\varepsilon$ for all $(t,j)\in\dom\phi$ with $t+j\geq T_1$. It means that the corresponding extended solution satisfies the uniform pre-attractivity conditions (\emph{ii}) of Definition~\ref{UGpASoZ} with $T=T_1$ and $K=0$.

It remains to show that the conditions (\emph{ii}) of Definition~\ref{UGpASoZ} are also satisfied for solutions starting outside the set $\A_1$. Let $|\tilde\phi(0,0,0)|_{\A_2}\leq r$, $\tilde\phi(0,0,0)\in\A_0\setminus\A_1$ and let the arc $\phi=\tilde\phi(t,j,0)$, $(t,j,0)\in\dom\tilde\phi$ be Zeno. From the conditions of Theorem \ref{main} its $\omega$-limit set is non-empty, hence this solution is being prolonged from the set $\Omega(\phi)$. From Lemma \ref{olimit} it follows that $\Omega(\phi)\subset \A_1$ and from Definition~\ref{phase} it follows that the set $\A_1$ can be considered as a new state space for the system $\H\cap\A_1$. Since $\A_2$ is UGpAS for the system $\H\cap\A_1$ and $\Omega(\phi)\subset \A_1$ it follows that an extended solution $\tilde\phi$ issued from $\A_0\setminus \A_1$ satisfies the pre-attractivity conditions (\emph{ii}) of Definition \ref{UGpASoZ} with $T=T_1$ an $K=1$.

Since there are no other types of solutions to the system $\H$ starting outside $\A_1$, the set $\A_2$ is UGpAoZ with $T=T_1$, $K=1$. UGpASoZ follows from UGSoZ and UGpAoZ. Iterating the previous reasoning one can prove UGpASoZ for any finite $n$. This concludes the proof.
\end{pf}

The proven result is applicable only to a system with the arcs, issued outside the set $\A_i$, $i=1,\ldots,n-1$, that are Zeno. However it is easy to prove the stability and pre-attractivity for the case of non-Zeno hybrid arcs that reach the corresponding set $\A_i$ in a finite time (e.g. when $t+j$ is bounded). In this case we lose the uniformity of the pre-attractivity which means that the constants $T$ and $K$ that describe the time needed for a solution to reach the $\varepsilon$-neighbourhood of the set $\A_n$ depend on a particular solution.

\begin{thm}\label{main2}
Let there exist a finite sequence $\A_n\subset \A_{n-1}\subset\ldots\subset \A_1\subset \A_0= C\cup D$ such that $\A_i$ is UGpAS for the system $\H\cap\A_{i-1}$, $i=1,\ldots,n$ and for all initial values $\phi(0,0)\in \A_{i-1} \setminus \A_i$, $i=1,\ldots,n-1$ the corresponding solutions $\phi$ are either Zeno with non-empty $\omega$-limit sets or $\rge\phi\cap\A_i \not = \emptyset$. Then $\A_n$ is UGSoZ and GpAoZ.
\end{thm}
\begin{pf}
The proof repeats the reasonings of the Theorem~\ref{main}. UGSoZ proof is similar. However, we need to consider the second type of arcs (non-Zeno) in order to prove the pre-attractivity. Let $\phi$ be the solution issued from a point outside the set $\A_1$ such that $\rge\phi\cap \A_1\not = \emptyset$. It means that there exists $(t^*,j^*)\in\dom\phi$ such that $\phi(t^*,j^*)=\xi^*\in \A_1$. The further evolution of solution $\phi$ coincides with the evolution of the solution $\phi^*$ to the system $\H\cap\A_1$ with initial condition $\phi^*(0,0)=\xi^*$. Since the set $\A_2$ is UGpAS for the system $\H\cap\A_1$ it follows that the solution issued from $\xi\in\A_0\setminus \A_1$ satisfies the conditions (\emph{iii}) from Definition~\ref{UGpASoZ} with $T=t^*+j^*+T_1$, $K=0$. Note that $t^*$ and $j^*$ can be different for every solution $\phi$. Denote $T_2(\varphi)=t^*+j^*$.

Since there are no other types of solutions to the system $\H$ issued outside $\A_1$, the set $\A_2$ is GpAoZ with $T=T_1+T_2(\phi)$, $K=1$. Iterating the previous reasoning one can prove GpAoZ for any finite $n$. This concludes the proof.
\end{pf}

Note that the situation when every set $\A_i$ is UGpAS but there exists an arc that wounds by spiral and tends to the set $\A_i$ but does not intersect it, does not satisfy the conditions of Theorem~\ref{main2}.

Next we present an example that demonstrates the usage of Theorem~\ref{main}. To check the UGpAS of a set $\A$ we will use the known theorem from \cite{goebel2012hybrid}:
\begin{proposition}[\cite{goebel2012hybrid}]\label{thTeel}
Let $\A\subset\R^n$ be closed. If $V(x)$ is a Lyapunov function candidate for $\H$ and there exist $\alpha_1, \alpha_2 \in \K_{\infty}$, and a continuous $\rho\in\PD$ such that
\begin{equation*}
\begin{array}{rll}
\alpha_1(|x|_{\A})\leq V(x)&\leq\alpha_2(|x|_{\A}) & \forall x\in C\cup D\cup g(D) \\
\langle \nabla V(x),f(x) \rangle &\leq -\rho(|x|_{\A}) & \forall x\in C \\
V(g(x))-V(x) &\leq -\rho(|x|_{\A}) & \forall x\in D
\end{array}
\end{equation*}
then $\A$ is UGpAS.
\end{proposition}

\begin{example}
Consider the following system with state $x=(x_1,x_2,x_3)\in\R^3$
\begin{equation}\label{ex3}
\begin{array}{lclllcll}
\dot x_1 & = & x_2, & ~ & {x_1}^+ & = & x_1, & ~\\
\dot x_2 & = & -\gamma(x_1,x_2), & x\in C, & {x_2}^+ & = & -\lambda x_2, & x\in D\\
\dot x_3 & = & -x_3, & ~ & {x_3}^+ & = & -x_3, & ~
%{x_1}^+ & = & x_1, & ~\\
%{x_2}^+ & = & -\lambda x_2, & x\in D\\
%{x_3}^+ & = & -x_3, & ~
\end{array}
\end{equation}
and with flow and jumps sets given by
\begin{equation*}
\begin{split}
C&=\{x\in\R^3: x_1>0 \text{~or~} x_1=0, x_2\geq0\}, \\
D&=\{x\in\R^3: x_1=0, x_2<0\},
\end{split}
\end{equation*}
where $\lambda\in(0,1)$ and $\gamma:\R^2\to\R$ is given by \eqref{eqgamma}.
\end{example}
Let us prove that the origin is UGpASoZ. As one may notice, flow and jump sets of system \eqref{ex3} do not depend on $x_3$. This system can be interpreted as an interconnection of a bouncing ball with state $(x_1,x_2)$ and some other process with state $x_3$. Each time when the ball bounces at the floor the state variable $x_3$ changes its sign.

Let function $V$ be defined by
\begin{equation*}%\label{V}
V(x)= (1+\theta\arctan{x_2})\left(\frac{x_2^2}{2}+\gamma x_1\right)
\end{equation*}
with
\begin{equation*}
\theta=\frac{1-\lambda^2}{\pi(1+\lambda^2)}\quad\text{and}\quad \gamma=9,81.
\end{equation*}

The set $\A_1=\{(x\in\R^3:x_1=x_2=0)\}$ is UGpAS since $V$ satisfies Proposition~\ref{thTeel} with respect to the origin for a single bouncing ball~\cite{goebel2012hybrid} and the distance from a point $(x_1,x_2,x_3)\in\mathbb R^3$ to the set $\A_1$ coincides with the Euclidean norm $\left\|\cdot\right\|$ of the corresponding vector $(x_1,x_2)\in\mathbb R^2$:$$|(x_1,x_2,x_3)|_{\A_1}=\left\|(x_1,x_2)\right\|=\sqrt{x_1^2+x_2^2}.$$

Then we arrive to a system of the form \eqref{ex3} with the new state space $\A_1$, the new flow set $C\cap\A_1=\{x\in\A_1:x_1>0 \cup x_1=0, x_2\geq0\}=\A_1$ and the new jump set $D\cap\A_1:=\{x\in\A_1: x_1=0, x_2<0\}=\emptyset$. This system is purely continuous and the origin $(0,0,0)$ is UGpAS. Since all the arcs issued outside the set $\A_1$ are Zeno, from Theorem~\ref{main} it follows that the origin is UGpASoZ.$\hfill\square$

Theorem~\ref{main} proposes a sequential narrowing of the state space of a hybrid system. For the last example this process can be described with the following sequence of sets:
\begin{figure}[h]
  \center
  \includegraphics[width=0.5\textwidth]{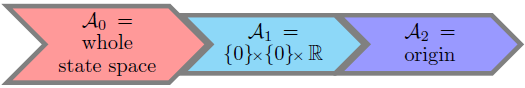}
\end{figure}

Theorem~\ref{main} has been used to prove UGpASoZ of the origin without constructing solutions to the system \eqref{ex3} and their prolongation explicitly. One may check that for a particular initial data the corresponding solutions to the system \eqref{ex3} have two $\omega$-limit points. If $\phi(0,0,0)=(a,b,c)\in\R^3$ with $\sqrt{a^2+b^2}\not =0$, $c\not =0$ then system \eqref{ex3} has Zeno hybrid arc with Zeno time $\tau=\frac{1}{g}\left(b+\frac{2\lambda}{1-\lambda}\sqrt{b^2+2ga} \right)$ and the $\omega$-limit set consisting of two points $(0,0,\pm c\cdot e^{-\tau})$. Hence the initial point $(a,b,c)$ generates two solutions. Despite such complex situation we were able to use Theorem~\ref{main} to verify UGpASoZ without knowing the exact number of $\omega$-limit points of Zeno arcs. Moreover, Theorem~\ref{main} along with Lemma~\ref{olimit} can be used for $\omega$-limit points localization. If one finds a function $V$ satisfying Proposition~\ref{thTeel} for some set $\A$, then following Lemma~\ref{olimit}, $\omega$-limit points of solutions are contained in the set $\A$.

\section{Discussion and open questions}
The results presented here are beneficial for construction and stability analysis of solutions to hybrid dynamical systems that exhibit Zeno behavior. The main contributions of the paper are the following. First, we have introduced the extended hybrid time domain and new prolonged solution concept that heavily relies on the axiomatics and notation of \cite{goebel2012hybrid}. These extended solutions helped us to avoid such undesired effects as freezing of solutions.
Second, we propose a generalisation of the attractivity concept and prove theorems that provide Lyapunov-like sufficient conditions for stability without knowing the explicit solution and $\omega$-limit points. However, in order to apply these theorems one should be able to verify whether all hybrid arcs are either Zeno or intersect an appropriate set $\A_i$. %Third, the up to date literature overview in Section~\ref{ss1} provides the reader with current state of the art in the field of hybrid systems exhibiting Zeno solutions and problems related to such behavior.

We believe that the proposed way of solution's prolongation from its $\omega$-limit points can also be achieved without the introduction of a new $3$-dimensional hybrid time domain. However it would cause a significant redefining of the basic concepts of hybrid dynamical systems framework. An important advantage of the proposed approach is the ability to utilize a wide range of previously developed results on UGpAS, e.g.~from~\cite{goebel2012hybrid}, for stability analysis of extended solutions.

Several problems have no answers yet and are very exciting to be solved. The first one is an extension of the results to infinite dimensional setting. This can be described by a vacuous interconnection of infinitely many bouncing balls. One can easily check that this system has a qualitatively different behavior depending on the initial conditions. Let the balls are enumerated by index $n$ from $1$ to $\infty$ and each ball starts its way with zero velocity and vertical position equals to $n$. Then the Zeno index $k$ of the hybrid time domain for this case tends to infinity while the ordinary time $t\to\infty$. However if every ball starts with the position $\frac{1}{2^n}$ then the Zeno index $k$ reaches infinity by a finite ordinary time $t$. If we interconnect each of these systems with a purely continuous process that tends to zero (like $\frac{dx}{dt}=-x$) then a solution of the entire interconnection will tend to the origin in the first case and will "freeze" away from the origin in the second one. This situation gives an intuition that such kind of systems can be treated using some analogues of a local stability concept and needs a deeper investigation for a comprehensive analysis of its behavior.

Another challenging issue is an interconnection of a completely continuous and a completely discrete system. The resulting flow and jump sets obtained in ''natural'' manner lead to a system with only discrete time domain. The examples of such processes are for instance sample-and-hold control where a discrete-time algorithm measures the state of a continuous time system and updates it. In this case an entire interconnected system will have a solution with only discrete time and we just lose the continuous process.

%\begin{ack}
%This work was supported by the German Federal Ministry of Education and Research (BMBF) as a part of the research project ''LadeRamProdukt''.
%\begin{figure}[hb]
%  \includegraphics[width=1.5in]{bild.jpg}
%  \hfill
%  \includegraphics[width=1.5in]{BMBF-FH.jpg}
%\end{figure}
%\end{ack}

%\vspace{-1mm}

\bibliography{mybibfile}             % bib file to produce the bibliography
                                                     % with bibtex (preferred)

\end{document}